\documentclass[11pt]{article}
\usepackage{amsmath,amsthm}

\def\Llr{\Longleftrightarrow}
\def\refeq#1{\if\workingver y(\ref{#1})-[[#1]]\else(\ref{#1})\fi}
\def\refth#1{\if\workingver y\ref{#1}-[[#1]]\else\ref{#1}\fi}
\def\mylabel#1{\if\workingver y\label{#1}{\bf\ \ [[#1]]\ \ }
\else\label{#1}\fi}
\def\mybibitem#1{\if\workingver y\bibitem{#1}{\bf\ \
[[#1]]\ \ }\else\bibitem{#1}\fi}

\newtheorem{thm}{Theorem}
\newtheorem{lem}[thm]{Lemma}
\newtheorem{defn}[thm]{Definition}
\newtheorem{cor}[thm]{Corollary}
\newtheorem{prop}[thm]{Proposition}

\newtheorem{rem}[thm]{Remark}

\let\workingver=n

      \advance\textheight\topskip
\font\tenrm=cmr10
\font\tensl=cmsl10
\font\ninerm=cmr9
\font\erm=cmr8

\def\eps{\epsilon}
\def\ga{\gamma}
\def\be{\beta}
\newcommand{\bb}{\bar}
\newcommand{\hs }{\hspace{.6cm}}
\newcommand{\PS}{\centerline{\scshape Pantelimon St\u anic\u a}}
\newcommand{\beq}{\begin{equation}}
\newcommand{\gata}{\end{equation}}
\def\st{\star}
\def\Mod{{\rm Mod}}
\def\dbl{\left[\left[}
\def\dbr{\right]\right]}
\def\Det{{\rm Det}}
\def\ZZ{V_n}
\def\Z{V_1}
\def\Zz{{\bf Z}_2^{n-1}}
\def\op{\oplus}
\def\b{\bar }
\def\co{{\cal O}}

\begin{document}
\title{Nonlinearity, Local and Global Avalanche Characteristics
of Balanced Boolean Functions}
\author{Pantelimon St\u anic\u a \vspace{.5cm}
\thanks{On leave from the Institute of
Mathematics of Romanian Academy,
Bucharest, Romania}\\
\em
\small Auburn University Montgomery, Department of Mathematics\\
\em
\small Montgomery, AL 36117, {e-mail: stanpan@strudel.aum.edu.}}
\date{}
\maketitle

\hrule

\vspace{.3cm}
{\noindent\bf Abstract}

\vspace{.2cm}
{\tenrm
For a Boolean function $f$, define
$\Delta_f(\alpha)=\sum_{{x}}{\hat f ({ x})}
{\hat f (x\op \alpha)}$, $\hat f(x)=(-1)^{f(x)}$,
the {\em absolute} indicator
\(
\Delta_f=\max_{\alpha\not= 0} |\Delta_f(\alpha)|,
\)
and  the {\em sum-of-squares} indicator
\(
\sigma_f=\sum_{\alpha} \Delta_f^2(\alpha).
\)
We construct a class of functions with good local
avalanche characteristics,
but bad global avalanche characteristics, namely
we show that $2^{2n}(1+p)\leq \sigma_f\leq 2^{3n-2},
\Delta_f=2^n$, where $p$ is
the number of linear structures (with even Hamming weight)
of the first half of an {\em SAC\/} balanced Boolean function
$f$. We also derive some bounds for the nonlinearity of
such functions. It improves upon the
results of Son {\em et al.} \cite{SLCS}
and Sung {\em et al.} \cite{SCP}. In our second result
we construct a class of highly nonlinear balanced
functions with good local and
global avalanche characteristics. We show that for these functions,
$2^{2n+2}\leq \sigma_f\leq 2^{2n+2+\epsilon}$
($\epsilon=0$ for $n$ even
 and $\epsilon=1$ for $n$ odd).
}

\vspace{.4cm}

{\ninerm
}

\vspace{.3cm}
\hrule
\pagestyle{myheadings}

\vspace{1.4cm}

\baselineskip=1.95\baselineskip

\section{Definitions and Preliminaries}

The design and evaluation of cryptographic functions
requires the definition
of design criteria.
The Strict Avalanche Criterion ({\em SAC\/}) was
introduced by Webster and
Tavares \cite{WT} in a study of these criteria.
A Boolean function is said to
satisfy the {\em SAC} if complementing a single bit results
in changing the output bit with probability exactly one
half. In \cite{PLLGV},
Preneel {\em et al.} introduced
the {\em propagation criterion of degree} $k$ ($PC$ of
degree $k$ or $PC(k)$),
which generalizes the {\em SAC}: a function
satisfies the $PC(k)$ if by complementing at most $k$
bits the output changes
with probability
exactly one half.  Obviously $PC(1)$ is equivalent to
the {\em SAC} property.
The $PC(k)$ can be stated in terms of autocorrelation
function.
Let $V_n=\{\alpha_i|1\leq i\leq 2^n\}$ be the set of vectors
of ${\bf Z}_2^n$ in lexicographical order.
 For a function
on $\ZZ$, we say
that $f$ satisfies the $PC(k)$ if and only if
\beq
\label{sac}
\sum_{x\in \ZZ} f(x)\oplus f(x\op c)=2^{n-1},
\end{equation}
for all elements $c$ with {\em Hamming weight} (the number of
nonzero bits) $1\leq wt(c)\leq k$,
or equivalently,  $\Delta_f(c)=0$, where
\[
\Delta_f(c)=\sum_{{x}\in \ZZ}{\hat f ({ x})}
{\hat f (x\op c)}
\]
is the autocorrelation function and $\hat f(x)=(-1)^{f(x)}$.
There is also another variation of the {\em PC}, when one
 requires to have
the  above relation for an arbitrary subset of $\ZZ$, not
necessarily for {\em all}
$x$ with $1\leq wt(x)\leq k$ (see also \cite{F}).

As many authors observed, the {\em PC} is a very important
concept in designing cryptographic primitives
used in data encryption algorithms and hash functions.
However, the {\em PC}
captures only local properties
of the function. In order to improve the global analysis
of cryptographically
strong functions, Zhang and Zheng \cite{ZZ1} introduced
another criterion, which
measures the {\em Global Avalanche Characteristics}
({\em GAC}) of a Boolean
function. They proposed two indicators related
to the {\em GAC}: the {\em absolute} indicator
\[
\Delta_f=\max_{\alpha\not= 0} |\Delta_f(\alpha)|,
\]
and  the {\em sum-of-squares} indicator
\[
\sigma_f=\sum_{\alpha} \Delta_f^2(\alpha).
\]

The smaller $\sigma_f,\Delta_f$ the better the {\em GAC}
of a function.
Zhang and Zheng
obtained some bounds on the two indicators:
\[
2^{2n}\leq \sigma_f\leq 2^{3n}, 0\leq \Delta_f\leq 2^n.
\]
The upper bound for $\sigma_f$ holds if and only if $f$ is
affine and the
lower bound holds if and only if $f$ is {\em bent} (satisfies the PC
with respect to all $x\not= 0$).

There is an interest in computing bounds of the two
indicators for various
classes of Boolean
functions.
Recently, Son, Lim, Chee and Sung \cite{SLCS} proved
\beq
\label{firstbound}
\sigma_f\geq 2^{2n}+2^{n+3},
\end{equation}
when $f$ is a balanced Boolean function, and Sung, Chee and Park
\cite{SCP}
proved that
if $f$ also satisfies the {\em PC} with respect to
 $A\subset \ZZ$, $t=|A|$,
then
\beq
\label{secondbound}
\sigma_f \geq{
\begin{cases}
&2^{2n}+2^6(2^n-t-1), \ \text{if}\ 0\leq t\leq 2^n-2^{n-3}-1,
t\ \text{odd}\\
&2^{2n}+2^6(2^n-t+2), \ \text{if}\ 0\leq t\leq 2^n-2^{n-3}-1,
t\ \text{even}\\
&\Big(1+\frac{1}{2^n-1-t} \Big)2^{2n},\ \text{if}\
2^n-2^{n-3}-1< t\leq 2^n-2.
\end{cases}
}
\end{equation}
The result (\ref{secondbound}) improves upon (\ref{firstbound}).
Using the above result the authors of \cite{SCP}
have derived some new bounds
for the nonlinearity
of a balanced Boolean function satisfying the
{\em PC} with respect to $t$ vectors.
We will improve their results significantly.

\newpage

We need the following
\begin{defn}
\item[1.]
 We call $e_i$ the $i$-th {\em basis vector} of $V_n$.
\item[2.]
 An {\em affine} function is a Boolean function of the form
$\displaystyle f(x)=\displaystyle \oplus_{i=1}^n c_i x_i\oplus c.$
$f$ is called {\em linear} if $c=0$.
\item[3.]
 The {\em truth table} of $f$ is the binary sequence
$
f=(v_1,v_2,\ldots,v_{2^n}),$ where $v_i=f(\alpha_i)$.
\item[4.]
 The {\em Hamming
weight} of a binary vector $v$, denoted by $wt(v)$ is defined as the
number of ones it contains. The {\em Hamming distance} between two
functions $f,g:V_n\to V_1$, denoted by $d(f,g)$ is defined as
$wt(f\oplus g)$. $f$ is {\em balanced\/} is $wt(f)=2^{n-1}$.
\item[5.]
 The {\em nonlinearity} of a function $f$, denoted by
$N_f$ is defined as
\(
\displaystyle
\min_{l\in A_n} d(f,l),
\)
where $A_n$ is the class of all affine function on $V_n$.
\item[6.]
A vector $0\not=\alpha\in V_n$ is a {\em linear structure} of $f$
if $f(x)\oplus f(x\oplus\alpha)$ is constant for all $x$.
\item[7.]
If $X,Y$ are two strings of the same length, $(X|Y)$
means that $X$ and $Y$
occupy the same positions in the first and the second
half of some function.
\item[8.]
Define the set of $4$-bit blocks
\(
T=\{
A=0,0,1,1;\ {\b A}=1,1,0,0;\ B=0,1,0,1;\ {\b B}=1,0,1,0;\
C=0,1,1,0;\ {\b C}=1,0,0,1;\ D=0,0,0,0;\ {\b D}=1,1,1,1\}.
\)
\item[9.] If some bits of an affine function $l$ agree with the
the corresponding bits in a function $f$, we say that $l$
{\em cancels} those bits in $f$.
\item[10.]
If $u$ is a given string and $g$ is a
 Boolean function, we use
$u^g=$ {\em the string of bits in $g$ which occupy the same positions as
the bits in the string $u$}.
\item[11.]
If a Boolean string is a concatenation of either $A/\b A$ or
$B/\b B$ or $C/\b C$ or $D/\b D$ we say that it is
 {\em based on} $A$ or $B$ or $C$ or $D$.
 \item[12.]
 By $MSB(\cdot)$ we denote the {\em most significant bit} of
the enclosed argument.
\end{defn}

\section{The First Result}

In this section the function $f$ will denote a
balanced Boolean function which
satisfies the {\em SAC}.
We will consider {\em SAC} functions constructed using some ideas of
\cite{YCST,YT} (see also
\cite{CS} for another version of the construction).
Define $\displaystyle {\bf 1}\cdot x=\oplus_{i=1}^{n-1} x_i$, if
$x=(x_1,\ldots,x_{n-1})$.
Let  $g:V_{n-1}\to V_1$ denote the Boolean function
${\bf 1}\cdot x\op b,\ b\in V_1$,
which satisfies
$g(x)=\bar g(x\op a)$,
for any element $a$ of odd Hamming weight.
For a vector $v\in V_n$, we denote by
$v'\in V_{n-1}$ the $n-1$ least
significant bits in $v$.
In \cite{YCST,YT,CS} or \cite{Stanica-Ph.D.}
it is proved that functions of the form
\beq
\label{concat1}
f=(h\ |h\op g),\ \text{or}\ f=(h\ |l\op g),
\end{equation}
are {\em SAC} functions, where $h$ is an arbitrary
 function on $V_{n-1}$
and
$l(x)=h(x\op a)$, $wt(a)=odd$. Let $\b x$ be the complement of $x$.

\begin{prop}
The functions \refeq{concat1} can be written as
$f(x_1,\ldots,x_{n-1},x_n)=$
\[
\begin{split}
& \b x_n h(x_1,\ldots,x_{n-1})\oplus  x_n \left(h(x_1,\ldots,x_{n-1})
\oplus_{i=1}^{n-1} x_i\oplus b\right)\ \text{or}\\
& \b x_n h(x_1,\ldots,x_{n-1})\oplus x_n
\left(h(x_1,\ldots,\b x_k,\ldots, x_{n-1})
\oplus_{i=1}^{n-1} x_i\oplus b\right),
\end{split}
\]
(an odd number of input bits $x_k$ are complemented), for an arbitrary
Boolean function $h$ defined on $V_{n-1}$ and $b\in V_1$.
\end{prop}
\begin{proof}
Straightforward using the definition of $g$ and concatenation.
\end{proof}

First, we consider the case of balanced Boolean functions
$f$ defined on $V_n, n\geq 3$ of the form (\ref{concat1})
such that $h$ has linear structures.  We denote by
${\cal L}_h^{even}$ the
number of
nonzero  linear structures of $h$ with even Hamming weight.
We take $a$ to be an element of odd Hamming weight.
In our next theorem we compute the indicators for
 a class of functions satisfying the SAC. We remark that the global
 characteristics are not good for these functions
  although the local ones are (the functions are SAC).

\begin{thm}
\label{thm}
If $f$ is a balanced Boolean function of the form $f=(h|l\op g)$,
$l(x)=h(x)$
or $l(x)=h(x\op a)$,
$h$ an arbitrary Boolean function with ${\cal L}_h^{even}\geq 1$
and $g$ as before, we have
\beq
\label{bound3}
2^{2n}(1+{\cal L}_h^{even})\leq \sigma_f\leq 2^{3n-2}.
\end{equation}
\end{thm}

\begin{proof}

Zhang and Zheng \cite{ZZ2} proved that for functions satisfying the
{\em SAC}, the nonlinearity satisfies
\beq
\label{2n2}
N_f\geq 2^{n-2}.
\end{equation}
In \cite{SLCS} the following inequality is obtained:
\beq
\label{nonlin_sigma}
N_f\leq 2^{n-1}-\frac{1}{2}\sqrt{\sigma_f/2^n}.
\end{equation}
Using (\ref{2n2}) and (\ref{nonlin_sigma}) we obtain easily the right
inequality of (\ref{bound3}),
that is
\[
\sigma_f\leq 2^{3n-2}.
\]

From the proof of Lemma 1 of \cite{SCP} we get that $\sigma_f$
satisfies
\[
\sigma_f=\sum_x \Delta_f^2(x)=2^6\sum_{x} (b_x-2^{n-3})^2+2^{n+4}
\sum_x (b_x-2^{n-3}),
\]
where $b_x=\frac{1}{2}\sum_y f(y)f(y\op x)$.
Using the trivial identity $a b=\frac{1}{2}(a+ b-a\op b)$
and the fact that $f$
is balanced, we get
\(
b_x=\frac{1}{4} \sum_y \left(
f(y)+f(y\op x)-f(y)\op f(y\op x)
 \right)=2^{n-2}-\frac{1}{4}\sum_y f(y)\op f(y\op x).
\)
We note that $f$ satisfies the {\em PC} with respect to $x$
if and only if
$b_x=2^{n-3}$.
Since $f$ is balanced, $\sum_x (b_x-2^{n-3})=0$.
It follows that
\[
\sigma_f=2^{2n}+2^6\sum_{wt(x)\geq 2}(b_x-2^{n-3})^2.
\]

We want to evaluate $\sum_{wt(x)\geq 2}(b_x-2^{n-3})^2$.
In order to do that we have to compute
\[
S_{x}=\sum_{y\in V_n} f(y)\op f(y\op x).
\]
{\em Case} 1: $MSB(x)=0$.\\
In this case
\beq
\label{s_x}
\begin{split}
S_{x}
=&\sum_{y\in \ZZ} f(y)\op f(y\op x)=\sum_{i=1}^{2^{n-1}}
h(v_i')\op h(v_i'\op x')+\\
&\sum_{i=1}^{2^{n-1}} h(v_i')\op h(v_i'\op x')\op g(v_i')\op
g(v_i'\op x').
\end{split}
\end{equation}
{\em Case} 1.1: $wt(x')=even$.\\
In this case, since $g$ satisfies $g(x)=\bar g(x\op a)$
for any element
with odd Hamming weight, it follows that $g(v_i'\op x')=g(v_i')$.
 Therefore,
the equation (\ref{s_x}) becomes
\[
S_{x}=2\sum_{i=1}^{2^{n-1}} h(v_i')\op h(v_i'\op x').
\]
When $x'$ is a linear structure of $h$, $S_x=2^n c$, where
$c=h(0)\op h(0\op x')$.\\
{\em Case} 1.2: $wt(x')=odd$.\\
Then $g(v_i'\op x')=\bar g(v_i')$ and (\ref{s_x}) becomes
\[
S_{x}=\sum_{i=1}^{2^{n-1}} h(v_i')\op h(v_i'\op x') +
\sum_{i=1}^{2^{n-1}} h(v_i')\op h(v_i'\op x')\op 1=2^{n-1}.
\]
{\em Case} 2: $MSB(x)=1$.\\
In this case, $S_x$ can be evaluated as follows:
\[
S_{x}=\sum_{i=1}^{2^{n-1}} h(v_i')\op h(v_i'\op x')
\op g(v_i'\op x')+
\sum_{i=1}^{2^{n-1}} h(v_i')\op h(v_i'\op x')\op g(v_i').
\]
{\em Case} 2.1: $wt(x')=even$.\\
Since $g(v_i')=g(v_i'\op x')$, we get
\[
S_x=2\sum_{i=1}^{2^{n-1}} h(v_i')\op h(v_i'\op x')\op g(v_i').
\]
{\em Case} 2.2:  $wt(x')=odd$.\\
Since $g(v_i'\op x')=\bar g(v_i')$, we get
\[
S_{x}=\sum_{i=1}^{2^{n-1}} h(v_i')\op h(v_i'\op x') +
\sum_{i=1}^{2^{n-1}} h(v_i')\op h(v_i'\op x')\op 1=2^{n-1}.
\]

From the above analysis we deduce that:\\
Case 1.1: $b_x=2^{n-2}-2^{-2}S_x,$ and if $x'$ is a
linear structure for $h$,
$b_x=2^{n-2}$ or $b_x=0$.\\
Case 1.2: $b_x=2^{n-3}$.\\
Case 2.1: $b_x=2^{n-2}-2^{-2}S_x$, and if $x'$ is a
linear structure for $h$,
$b_x=2^{n-3}.$\\
Case 2.2: $b_x=2^{n-3}$.

We observe that the only cases where we do not know
precisely $b_x$ are
when $x$ is an element of odd Hamming weight with $x'$ not a linear
structure for $h$.

We deduce that in the case 1.1  with $x'$ a linear structure for $h$,
\[
(b_{x}-2^{n-3})^2=2^{2(n-3)}.
\]

Now, returning to the computation of $\sigma_f$,
with the new results we get
\[
\begin{split}
\sigma_f
=& 2^{2n}+2^6\sum_{wt(x)\geq 2} (b_x-2^{n-3})^2\geq \\
& 2^{2n}+2^6  2^{2(n-3)} {\cal L}_h^{even}=2^{2n}
\left(1+{\cal L}_h^{even}\right).
\end{split}
\]
\end{proof}

With the same data as in the previous theorem we obtain
\begin{cor}
For $n\geq 3, \Delta_f=2^n$.
\end{cor}

\begin{proof}
The corollary follows from the proof of the theorem.
For a Boolean balanced function, $\Delta_f(x)=2^3b_x-2^n$.
Therefore for any $x$,
such
that $x'$ is a linear structure of $h$ of even Hamming weight,
we have
\(b_x=0\ \text{or}\ 2^n\). Thus
$\Delta_f=\max_{x\in\ZZ} |\Delta_f(x)|=2^n.$
\end{proof}

The previous corollary can also be deduced from Lemma 7 of \cite{ZZ1},
observing that if $x'$ is a linear structure of $h$ with
even Hamming weight, then
$(0,x')$ is a linear structure for~ $f$.

The following is an easy consequence of the previous theorem.
It shows that the theorem
gives tight bounds.
\begin{cor}
\label{cor1}
For a balanced Boolean {\em SAC} function $f$ given by
$(\ref{concat1})$,
where $h$ is affine
we have the following equation
\[
\label{sigma-eq}
\sigma_f= 2^{3n-2}.
\]
\end{cor}

\begin{proof}
This follows from the fact that any nonzero element of $\ZZ$
is a linear
structure for an affine function.
\end{proof}

Now we turn our attention to the nonlinearity of such functions.
Using
\[
N_f\leq 2^{n-1}-2^{-n/2-1}\sqrt{\sigma_f},
\]
and $\sigma_f\geq 2^{2n}(1+{\cal L}_h^{even})$, we get the corollary

\begin{cor}
Let $f$ be as in the Theorem \ref{thm}. Then, the nonlinearity
satisfies
\beq
\label{nonlin_bound1}
2^{n-2}\leq N_f\leq 2^{n-1}-2^{n/2-1}\sqrt{1+{\cal L}_h^{even}}.
\end{equation}
If $f$ satisfies the conditions of Corollary \ref{cor1}, then we have
\beq
\label{nonlin_bound2}
N_f= 2^{n-2}.
\end{equation}
\end{cor}

Since
\(
2^n+2^{n/2+3}+2^4<2^n\left(1+{\cal L}_h^{even}\right),\ \text{if}\
{\cal L}_h^{even}\geq 1,
\)
it follows that the bounds (\ref{nonlin_bound1}) or
(\ref{nonlin_bound2})
are better than the result of Zhang and Zheng,
who proved in \cite{ZZ2} that
\[
N_f\leq 2^{n-1}-\frac{1}{2}\sqrt{2^n+2^{n/2+3}+2^4},\ \text{if}\ n\
\text{is even}.
\]

Sung {\em et al.} \cite{SCP} obtained the following upper bound
for the nonlinearity
\[
\begin{split}
& N_f\leq 2^{n-1}-\frac{1}{2}\sqrt{2^n+2^6-\frac{(n+1)2^{6}}{2^n}},\
\text{if $n>2$ is odd and}\\
& N_f\leq 2^{n-1}-\frac{1}{2}\sqrt{2^n+2^6 - \frac{(n-1)2^6}{2^n}},\
\text{if $n$ is even},
\end{split}
\]
which is certainly weaker than the bound we have obtained.


\section{Highly nonlinear balanced SAC functions with good GAC}

In the previous section we constructed a class of balanced
functions with good local
avalanche characteristics, but bad global avalanche characteristics.
In this section we will use some results from \cite{Stanica-Ph.D.} to
construct balanced
Boolean {\em SAC} functions of nonlinearity at least $2^n-2^{[(n+1)/2]}$,
with good {\em GAC}.

From a result we like to call {\em Folklore Lemma} (see
\cite{Stanica-Ph.D.}), we
know that
for any affine function $l$, if $L$ is the first string  of
length $2^s$ in
$l$, then the next string
of the same length will be $L$ or $\b L$. A consequence of this fact
 is that any affine function is made up as a concatenation
of blocks $A/\b A$ or $B/\b B$ or $C/\b C$ or $D/\b D$.

\begin{sloppypar}
Our next theorem was proven initially in a more general form.
However, its proof relied heavily on results available only in
\cite{Stanica-Ph.D.}, so we decided to provide here a complete
proof for a slightly restricted subclass. Moreover, for this subclass
we can provide better results, especially for even dimensions, which
makes it all worthwhile.
For the purpose of easy computation,
we define a transformation ${\cal O}(g)$ ("opposite") which maps an
affine function based on $M\in T$, into an affine function
based on the same block $M$, having the self-invertible property
$\co\left(\co(g)\right)=g$.
If $g=X_1 X_2\ldots X_{2^{n-2}}$, then
$\co(g)=Y_1Y_2\ldots Y_{2^{n-2}}$ is constructed by the
following {\em Algorithm},
supported by the  Folklore Lemma:\\
{\em
{ Step 1.} $Y_1=X_1$.\\
{ Step $i+2$.} For any $0\leq i\leq n-3$, if
$X_{2^i+1} \ldots X_{2^{i+1}}=X_{1} \ldots X_{2^{i}}$,
then $Y_{2^i+1} \ldots Y_{2^{i+1}}=\b Y_{1}\ldots\b Y_{2^{i}}$.
If $X_{2^i+1} \ldots X_{2^{i+1}}=\b X_{1} \ldots\b
X_{2^{i}}$, then $Y_{2^i+1} \ldots Y_{2^{i+1}}= Y_{1}\ldots
Y_{2^{i}}$.
}
\end{sloppypar}
\begin{rem}
The results will not change if we take the first block $Y_1=\b X_1$.
\end{rem}
By induction we can easily prove
\begin{lem}
 $\co(\b g)=\overline{\co(g)}$.
\end{lem}
The following theorem is a construction for balanced functions
of high nonlinearity with very good local and global avalanche
characteristics. Let $[x]$  ({\em the floor function})
to be the largest integer less
than or equal to $x$. For easy writing we let $h_i=\co(g_i)$.

\begin{thm}
\label{thm2}
For $n=2k\geq 8\ (\text{or}\ n=2k+1\geq 9)$ let $f$ to be the function
obtained by
concatenating $2^{k-1}$  segments
$T_i$. For each $1\leq i\leq 2^{k-2}$,
$T_i$ is of the form
\beq
\label{s_i}
(g_i h_i  g_i\b h_i|\b h_i g_i  h_i g_i)
\end{equation}
and the segment $T_{i+2^{k-2}}$ is of the form
\beq
\label{next_s_i}
(h_i\b g_i \b h_i\b g_i|\b g_i \b h_i  \b g_i h_i),
\end{equation}
respectively, where the functions $g_i$ are affine
functions on $\Z^{k-2}$ (or $\Z^{k-1}$).
Furthermore, we impose the following conditions:\\
$(i)$ Exactly a quarter of the functions $g_i$
 are based on each of the $4$-bit blocks $A,B,C,D$.\\
$(ii)$
For any  $1\leq i\not = j\leq 2^{k-2}$, the functions
$g_i\op g_j$ are balanced.

Then the function $f$ is balanced, satisfies the {\em SAC}, has the
nonlinearity $N_f\geq 2^{n-1}-2^{[\frac{n+1}{2}]}$
and the sum-of-squares indicator satisfies
\[
2^{2n+2}\leq  \sigma_f\leq 2^{2n+2+\epsilon},
\]
where $\epsilon=0,1$ if $n$ is even, respectively, odd.
\end{thm}

\begin{proof}
We will prove the theorem for the case of $n$ even, that is $n=2k$,
pointing out,
whenever necessary, the differences for the case of odd $n$.
The function $f$ can be written as
\beq
\label{goodGAC}
\begin{split}
\big(
&g_1 h_1 g_1 \b h_1\cdots g_{2^{k-2}} h_{2^{k-2}} g_{2^{k-2}}
\b h_{2^{k-2}}\quad
h_1 \b g_1\b h_1\b g_1\cdots h_{2^{k-2}}\b g_{2^{k-2}}
\b h_{2^{k-2}}\b g_{2^{k-2}}\\
&\b h_1g_1h_1g_1 \cdots \b h_{2^{k-2}}g_{2^{k-2}}h_{2^{k-2}}
g_{2^{k-2}}\quad
\b g_1\b h_1\b g_1 h_1\cdots\b g_{2^{k-2}}\b h_{2^{k-2}}
\b g_{2^{k-2}} h_{2^{k-2}}
\big).
\end{split}
\end{equation}
The fact that $f$ is balanced can be seen by pairing the
functions $g$ with
$\b g$ and $h$ with $\b h$
 in the two segments $T_i$ and $T_{i+2^{k-2}}$.
To show that $f$ satisfies the {\em SAC} we use
 some results of
Cusick and St\u anic\u a,
that is Lemma 1 or relation (8) of \cite{CS}, which says that a
function $f=(v_1,\ldots,v_{2^n})=X_1\cdots X_{2^{n-2}}$
satisfies the {\em SAC}\/
if and only if
\beq
\label{sac2}
\begin{split}
&({w_1}{ w_{2^{i-1}+1}}+{w_2}{w_{2^{i-1}+2}}+\cdots+{w_{2^{i-1}}}
{w_{2^i}})+\\
&({w_{2^i+1}}{w_{2^i+2^{i-1}+1}}+\cdots+{w_{2^i+2^{i-1}}}
{w_{2^{i+1}}})+\cdots +\\
&(w_{2^n-2^i+1}w_{2^n-2^{i-1}+1}+\cdots+w_{2^n-2^{i-1}}w_{2^n}) =0,
\end{split}
\end{equation}
for each $i=1,2,\ldots,n$, where $w_i=(-1)^{v_i},$ or equivalently
(if $i\geq 3$),
\beq
\label{sac_op}
(X_1\odot X_{2^{i-3}+1}+\cdots+X_{2^{i-3}}\odot X_{2^{i-2}})+\cdots=0,
\end{equation}
for each $i=3,4,\ldots,n$, where $M\odot N$ is equal to the
number of 0's minus the number of 1's in $M\op N$
If we associate the 4-bit blocks $\{ A,\b A\}\Llr \{ B,\b B\}$ and
$\{ C,\b C\}\Llr \{D,\b D\}$, we
see that, for $i\leq 2$, the relation (\ref{sac2}) holds.
Obviously, if $M\oplus N$ is balanced, then $M\odot N=0$. Thus,
in the sum
(\ref{sac_op}) the sum in each parenthesis is zero,
except perhaps the ones based entirely on $D,\b D$ (which are the
only unbalanced 4-bit blocks in $T$).
However, those terms will have an antidote in another parenthesis.
For instance, since $D\odot D=-D\odot \b D=4$,
$D\odot D$ will have the antidote $D\odot \b D$, according to the
form of our functions.

In order to compute the nonlinearity of $f$ we have counted the
bits at which
our function differ from any linear or affine function.
Intuitively, we need to prove that on
average an affine function
cannot cancel to many blocks in a segment.
 Precisely, we show that given any two segments $U_1,U_2$
in the same half of $f$, based on the same block $M\in T$, then
$wt(U_1U_2\op U_1^lU_2^l)\geq 2^{k-1}+2^k$, for any affine function
$l$ based on the same block $M$. This is shown easily using the
folklore lemma, and observing that on the positions
of $U_1U_2$, $l$ can have only the following
forms: $(LLLL LLLL|LLLL LLLL)$,
$(LLLL LLLL|\b L\b L\b L\b L\b L\b L\b L\b L)$,
$(LL\b L\b L LL\b L\b L|LL\b L\b L LL\b L\b L)$, etc.
Since all cases are treated similarly, we may assume that
$(U_1^lU_2^l)=(LLLL LLLL|LLLL LLLL)$
(recall the definition of $U^l$).
Without loss of generality we may assume that $U_1,U_2$ are in
 the first half of $f$ and
$U_1=(g_1 h_1 g_1 \b h_1) |\b h_1 g_1 h_1 g_1 )$,
$U_2=(g_2 h_2 g_2 \b h_2) |\b h_2 g_2 h_2 g_2 )$.
 Thus
\begin{eqnarray*}
 wt(U_1 U_2\oplus U_1^lU_2^l)
 &=&
 2 wt(g_1\op L)+ wt(h_1\op L)+  wt(\b h_1\op L)\\
&&+
2 wt(g_2\op L)+ wt(h_2\op L)+ wt(\b h_2\op L)\\
&&+
wt(\b h_1\op L)+ wt(h_1\op L)+ 2 wt(g_1\op L)\\
&&+
wt(\b h_2\op L)+ wt(h_2\op L)+ 2 wt(g_2\op L)\\
&=& 4 wt(g_1\op L)+4 wt(g_2\op L)+2^{k}\\
&\geq &
4 wt(g_1\op g_2)+2^{k}= 2^{k-1}+2^{k}.
\end{eqnarray*}
Here we used
$wt(a\op c)+wt(b\op c)\geq wt(a\op b)$,
 the fact that $g_i\op g_j$ is balanced and
$wt(a\op b)+wt(a\op\b b)=2^{k-2}$, if $a,b,c\in V_{k-2}$.
Next, we compute $wt(f\op l)$.
One may assume that $l$ is based on $A$. From the part
of $f$ that does not contain $A, \b A$
we get $3\cdot 2^{2k-3}=2^{2k-1}-2^{2k-3}$ units for the weight
(we recall that only a quarter of all blocks contain  $A,\b A$).
 We consider now the part of $f$ based on $A$.
Using the previous result, we deduce that in the worst case
(minimum weight), $l$ cancels completely at
 most four functions from each half,
and from the rest of the part of $f$ based on $A$,
half of the blocks are cancelled.
Since there are $2^{k}$ functions based on $A$ and we cancel $8$
functions, we gather that
there remain $2^{k}-8$ functions uncancelled. Since each
uncancelled function contributes $2^{k-3}$ units to the weight
(recall that if two affine functions $g,l$ are not
equal or complementary, their sum is balanced),
 we get
$2^{2k-3}-2^k$ units contributed to the weight by the part based
on $A$, so the nonlinearity is at least
$2^{2k-1}-2^{2k-3}+2^{2k-3}-2^k=2^{2k-1}-2^k$.
 In the odd case we get $N_f\geq 2^{2k-1}-2^{k+1}$
 (the lengths of the affine functions $g_i,h_i$ double, while
the number of segments remains the same), by a similar argument.

Now, since
\(
N_f\leq 2^{n-1}-2^{-\frac{n}{2}-1}\sqrt{\sigma_f}
\)
and from the above analysis
\(
N_f\geq 2^{n-1}-2^{[\frac{n+1}{2}]}
\)
we get
\[
2^{n-1}-2^{[\frac{n+1}{2}]}\leq 2^{n-1}-
2^{-\frac{n}{2}-1}\sqrt{\sigma_f},
\]
which will produce our right hand side inequality
\[
\sigma_f\leq 2^{2n+2}, \ \text{if $n$ is even, and }\
\sigma_f\leq 2^{2n+3},\
\text{if $n$ is odd}.
\]

In order to evaluate $S_x$ for suitably chosen $x$ we apply the same
technique as in the
proof of Theorem \refth{thm}.
For $x=e_i\op e_j, i<j$, let
\beq
\label{s_ij}
\begin{split}
S_x =
&\sum_{y\in\ZZ} f(y)\op f(y\op x)=
\sum_{s=1}^{2^n} f(v_s)\op f(v_s\op e_i\op e_j)=\\
& 2 [f(v_1)\op f(v_{2^{j-1}+2^{i-1}+1}) +\cdots+
f(v_{2^{i-1}})\op f(v_{2^{j-1}+2^i})+\\
& f(v_{2^{i-1}+1})\op f(v_{2^{j-1}+1})+\cdots+
f(v_{2^{i-1}+2^{i-1}})\op f(v_{2^{j-1}+2^{i-1}})]+\cdots.
\end{split}
\end{equation}
Using the form of our functions and taking
$x=e_{n-1}\op e_{n}$, we get
\[
S_{e_{n-1}\op e_{n}}
=  2\sum_{g_i,h_i} (g_i\op \b g_i+h_i\op \b h_i+
g_i\op \b g_i+\b h_i\op h_i)=2^n.
\]
Thus,
$(b_{e_{n-1}\op e_{n}}-2^{n-3})^2=2^{2n-6}$.

Now, we take $x=e_i\op e_j\op e_r, i<j<r$.
Thus, we get
\beq
\label{s_ijr}
\begin{split}
S_x=
&\sum_{y\in\ZZ} f(y)\op f(y\op x)=\\
&\sum_{s=1}^{2^n} f(v_s)\op f(v_s\op e_i\op e_j\op e_r)=\\
& 2 [f(v_1)\op f(v_{2^{r-1}+2^{j-1}+2^{i-1}+1}) +\cdots+\\
&f(v_{2^{i-1}})\op f(v_{2^{r-1}+2^{j-1}+2^i})+\\
& f(v_{2^{i-1}+1})\op f(v_{2^{r-1}+2^{j-1}+1})+\cdots+ \\
&f(v_{2^{i-1}+2^{i-1}})\op f(v_{2^{r-1}+2^{j-1}+2^{i-1}})]+\cdots.
\end{split}
\end{equation}
Now, taking $x=e_{k-1}\op e_k\op e_n$ and $n=2k$, we obtain
\[\begin{split}
S_{e_{k-1}\op e_k \op e_n}=2\bigg[
&\big(f(v_1)\op f(v_{2^{n-1}+2^{k-1}+2^{k-2}+1}) +\cdots+ \\
&f(v_{2^{k-2}})\op f(v_{2^{n-1}+2^k})\big)+\\
&\big( f(v_{2^{{k-2}}+1})\op f(v_{2^{n-1}+2^{k-1}+1})+\cdots+ \\
&f(v_{2^{k-2}+2^{k-2}})\op f(v_{2^{n-1}+2^{k-1}+2^{k-2}})\big)+\\
&\big( f(v_{2^{{k-1}}+1})\op f(v_{2^{n-1}+2^{k-2}+1})+\cdots+ \\
&f(v_{2^{k-1}+2^{k-2}})\op f(v_{2^{n-1}+2^{k-1}})\big)+\\
&\big( f(v_{2^{k-1}+2^{k-2}+1})\op f(v_{2^{n-1}+1})+\cdots+ \\
&f(v_{2^k})\op f(v_{2^{n-1}+2^{k-2}})\big)
\bigg]+\cdots
\end{split}
\]
for any function $f$. In particular, for the functions in our class,
we get
\begin{eqnarray}
S_{e_{k-1}\op e_k \op e_n}
&=& 2 \sum_{s=1}^{2^{k-2}}(g_s\op  g_s+h_s\op  h_s+g_s\op
g_s+\b h_s\op\b
h_s)\nonumber\\
&+& 2 \sum_{s=1}^{2^{k-2}}(h_s\op h_s+\b g_s\op \b g_s+
\b h_s\op \b h_s+\b
g_s\op\b g_s)=0.\nonumber
\end{eqnarray}
Similarly, $S_{e_{k-1}\op e_k\op e_{n-1}}=2^n$.
Thus, $b_{e_{k-1}\op e_k\op e_n}=2^{n-2}$
 and $b_{e_{k-1}\op e_k\op e_{n-1}}=0.$

In any of the three cases
$x=e_{n-1}\op e_n,e_{k-1}\op e_k\op e_{n-1},e_{k-1}\op e_k\op e_{n}$,
we have $(b_x-2^{n-3})^2=2^{2n-6}$.
 Thus,
\[
\sigma_f\geq 2^{2n}+2^6 2^{2n-6}+2^6 2^{2n-6}+2^6 2^{2n-6}=
2^{2n+2}.
\]
\end{proof}

\begin{cor}
For $f$ given by Theorem \refth{thm2}, we have
\(
\Delta_f= 2^n.
\)
\end{cor}

\begin{proof}
We know that $\Delta_f(x)=2^3b_x-2^n$.
Therefore,
\[
\Delta_f(e_{k-1}\op e_k\op e_n)=2^3 \cdot 2^{n-2}-2^n=2^n,
\]
and the result follows.
\end{proof}

\begin{cor}
If $n$ is even and $f$ is given as in Theorem \refth{thm2}, then
$\sigma_f=2^{2n+2}$, $N_f=2^{n-1}-2^{\frac{n}{2}}$, and
$f$ is PC with respect to all but four vectors. Moreover, the
three nonzero vectors, which do not satisfy the propagation criterion,
are linear structures for $f$.
\end{cor}
\begin{proof}
We proved that, if $n$ is even, then $\sigma_f=2^{2n+2}$.
If there is an $x$ not equal to the four displayed vectors in
the proof of Theorem \refth{thm2}, for which $f$ is not PC,
then $b_x\not=2^{n-3}$.
If so, then by the same argument we would get $\sigma_f>2^{2n+2}$,
 which
is not true. So $f$ is PC with respect to all but four vectors.
In \cite{ZZ3}, Zhang and Zheng proved that, if a function satisfies
the PC with respect to all but four vectors, then $n$ must be even,
the nonzero vectors, where the propagation criterion is not satisfied,
must be linear structures and $N_f=2^{n-1}-2^{n/2}$.
We have the result.
\end{proof}

As we can see the bounds are extremely good,
not too far from that of bent functions,
improving upon any known ones.
We suspect we can modify the construction to improve
the nonlinearity for the odd dimension as well, and we will pursue
this idea elsewhere.

\begin{sloppypar}
\begin{rem}
If the conditions imposed in Theorem \refth{thm2}
hold for $g_i$, they certainly hold
for $h_i=\co(g_i)$ as well.
\end{rem}
\end{sloppypar}

\section{Examples and Further Research}

An example of a function satisfying the conditions of
 Theorem \refth{thm2} with $h_i=\co(g_i)$, for $n=8$ is
\[
\begin{split}
& AAA\b ABBB\b B CCC\b C DDD\b D
A\b A\b A\b AB \b B\b B\b B C \b C\b C\b C D \b D\b D\b D \\
& \b AAAA\b BBB B\b CCC C\b DDD D
\b A\b A\b AA\b B \b B \b B B\b C \b C \b C C\b D \b D \b D D,
\end{split}
\]
which is balanced, SAC (actually, it is
PC with respect to all but
${\mathbf 0}, e_7\op e_8,e_3\op e_4\op e_8, e_3\op e_4\op e_7$),
has nonlinearity 112 and
the sum-of-squares indicator attains the upper bound,
$\sigma_f=262,144=2^{2\cdot 8+2}$. The
algebraic normal form is
$ x_1 + x_7+x_1 x_5+ x_1 x_6  + x_2 x_5 + x_2 x_6  + x_3 x_8  +
x_4 x_7+ x_4 x_8+x_5 x_6$.

We can define the transformation $\co$ using the same algorithm
starting with the first bit, rather than the first block,
so $\co(A)=B, \co(C)=D$, etc., obtaining
 a result similar to our Theorem \refth{thm2}. It seems that the
 algebraic degree increases for that class, but we were not able to
 prove that in its full generality.
An example of a function constructed using this idea, for $n=8$, is
\[
\begin{split}
& A B A \bar B B A B \bar A C D C \bar D D C D \bar C B
\bar A\bar B\bar A A \bar B
\bar A\bar B C \bar D\bar C\bar D D\bar C\b D\b C\\
& \bar B A B A \b A B A B \bar D C D C \b C D C D
\b A\b B\b A B \b B\b A\b B A\b D\b C\b D C\b C\b D\b C D.
\end{split}
\]
It turns out that the above function is balanced, has nonlinearity
 precisely 112,
it is SAC (in fact, it is PC with respect to 252 vectors),
the sum-of-squares indicator attains the upper bound,
$\sigma_f=262,144=2^{2\cdot 8+2}$. The algebraic normal form is
$x_1+x_7+x_1 x_5+x_1 x_6+x_1 x_7+x_1 x_8+
x_2 x_5+x_2 x_6+x_2 x_7+x_2 x_8+x_3 x_8+x_4 x_7+x_4 x_8+
x_5 x_6+x_6 x_7+x_6 x_8+x_2 x_3 x_7+x_2 x_3 x_8.$

Another venue of further research would be the construction of a
class of functions
with these good local and global avalanche characteristics and high
nonlinearity, using blocks in the complementary set of $T$,
namely
$T'=
\{
U=1,0,0,0;\ {\b U}=0,1,1,1;\ V=0,0,0,1;\ {\b V}=1,1,1,0;
X=0,1,0,0;\ {\b X}=1,0,1,1;\ Y=0,0,1,0;\ {\b Y}=1,1,0,1
\}$. Our experiments showed that this approach
seems to increase the algebraic degree of
the functions involved, but we were not able to find and control all
the mentioned cryptographic parameters, yet.

\noindent{\bf Acknowledgements.}
{
The author would like to thank the anonymous referees for their
 helpful comments, which improved
significantly the presentation of the paper.
}

\end{document}